\documentstyle[12pt]{article}
\begin{document}
\title{ A note on potentially $K_4-e$ graphical sequences
\thanks{ The Project Sponsored by NSF of Fujian and Fujian Provincial Training
Foundation   for "Bai-Quan-Wan Talents Engineering".}}
\author{{Chunhui Lai}\\
{\small Department of Mathematics}\\{ Zhangzhou Teachers College,
Zhangzhou} \\{\small Fujian 363000,
 P. R. of CHINA.}\\{\small e-mail: zjlaichu@public.zzptt.fj.cn}}
\date{}
\maketitle
\begin{center}
\begin{minipage}{120mm}
\vskip 0.1in
\begin{center}{\bf Abstract}\end{center}
 {A sequence $S$ is potentially $K_4-e$ graphical if it has
a realization containing a $K_4-e$ as a subgraph. Let
$\sigma(K_4-e, n)$ denote the  smallest degree sum such that every
$n$-term graphical sequence $S$ with $\sigma(S)\geq \sigma(K_4-e,
n)$ is potentially $K_4-e$ graphical. Gould, Jacobson, Lehel
raised the problem of  determining  the value of $\sigma (K_4-e,
n)$. In this paper, we prove that $\sigma (K_4-e, n)=2[(3n-1)/2]$
for $n\geq 7$, and $n=4,5,$ and  $\sigma(K_4-e, 6)= 20$.}
\end{minipage}
\end{center}
\vskip 0.3in
\baselineskip 14pt
\section * {1.Introduction}
If $S=(d_1,d_2,...,d_n)$ is a sequence of non-negative integers, then it
is called  graphical if there is a simple graph $G$ of order $n$, whose degree
sequence $(d(v_1 ),d(v_2 ),...,d(v_n ))$ is precisely $S$. If $G$ is such
a graph then $G$ is said to realize $S$ or be a realization of $S$. A
graphical sequence $S$ is potentially $H$ graphical if there is a
realization of $S$ containing $H$ as a subgarph, while $S$ is forcibly $H$
graphical if every realization of $S$ contains $H$ as a subgraph. Let
$\sigma(S)=d_1+d_2+... +d_n,$ $[x]$ be the largest integer less than or equal
to $x$, If $G$ and $G_1$ are graphs, then $G\cup G_1$ is the disjoint union
of $G$ and $G_1$. If $G = G_1$, we abbreviate $G\cup G_1$ as $2G$. Let $K_k$ be a complete
graph on $k$ vertices, $C_k$ be a cycle of length $k$.
\par
   Given a graph $H$, what is $ex(n,H)$, the maximum number of edges of a graph
with $n$ vertices not containing $H$ as a subgraph? This problem was proposed
for $H = C_4$ by Erdos [2] in 1938 and in general by Turan [9]. In the terms
of graphic sequences, the number $2ex(n,H)+2$  is  the minimum even  integer
$m$ such that every $n$-term graphical sequence $S$ with $\sigma (S)\geq m $
is forcibly $H$ graphical. Here we consider the following variant: determine
the  minimum even integer $m$ such that every $n$-term graphical sequence $S$
with $\sigma(S)\geq m$ is potentially $H$ graphical, We denote this minimum
$m$ by $\sigma (H,n)$. Erdos, Jacobson and  Lehel [1] show that
$\sigma(K_k,n)\geq (k-2)(2n-k+1)+2$; and conjecture that $\sigma(K_k,
n)=(k-2)(2n-k+1)+2$; They proved that if $S$ does not contain zero  terms, this
conjecture is true for $k=3, n\geq6$, Li and Song [6,7,8] proved that if
$S$ does   not contain zero terms, this conjecture is true for $k=4, n\geq 8$
and $k=5, n\geq 10$, and $\sigma (K_k, n) \leq 2n(k-2) +2$ for $ n\geq2k-1$.
Gould, Jacobson and Lehel [3] proved that this  conjecture is true for $ k=4,
n\geq 9$, if $n=8$ and $\sigma(S)\geq 28$, then either there is a
realization of $S$ containing $K_4$ or $S=(4^7,0^1)$(i.e. $S$ consists of $7$
integers $4$ and $1$  integer $0$); $\sigma (pK_2, n)=(p-1)(2n-2)+2 $ for $ p\geq
2;\sigma (C_4, n)=2[(3n-1)/2]$ for $n\geq 4, \sigma (C_4, n) \leq
\sigma(K_4-e , n) \leq \sigma(K_4, n)$; and raised the problem of determining
the value of  $\sigma(K_4-e, n)$. Lai[4,5] proved that $\sigma(C_{2m+1},
n)=m(2n-m-1)+2$, for $m\geq 2, n\geq 3m$; $\sigma(C_{2m+2} , n)=m(2n-m-1)+4$,
 for $ m\geq 2, n\geq 5m-2$. In this paper, we determine the values of
$\sigma(K_4-e, n)$.
 \section * {2.$\sigma(K_4-e,n)$}
 {\bf  Theorem 1.} For $n=4,5$ and $n\geq 7$
 \[\sigma(K_4-e,n)=\left\{ \begin{array}{cc}3n-1 & \mbox {if $n$
is odd}  \\3n-2 & \mbox {if $n$ is even.}\end{array}\right.\] For
$ n=6,$ $S$ is a 6-term graphical sequence with $\sigma(S) \geq
16$, then either there is a  realization of S containing $K_4-e$
or $S=(3^6)$. (Thus $\sigma(K_4-e, 6)=20)$.
\par
Proof. By [3],\[\sigma(K_4-e, n)\geq \sigma(C_4, n)= \left \{\begin{array}
{cc}3n-1 & \mbox {if $n$ is odd}\\3n-2 & \mbox { if $n$ is
even.}\end{array}\right.\]
For $ n\geq 4$. Assume $d_1 \geq d_2 \geq ...\geq d_n \geq 0$.
\par
For $n=4$, if a graph has size $q\geq 5$, then clearly it contains a $K_4-e$,
so that $\sigma (K_4-e,n)\leq 3n-2$.
\par
For $n=5$, we have $q\geq 7$. There are exactly $4$ graphs of order $5$ and size
$7$ and each  contains a $K_4-e$. Thus $\sigma (K_4-e,n)\leq 3n-1$.
\par
Suppose for $5\leq t<n, S_1$ is a $t$-term graphical sequence such that
 \[\sigma (S_1)\geq \left \{\begin{array}{cc}3t-1 & \mbox {if $t$ is odd}
 \\3t-2 & \mbox{if $t$ is even.}\end{array}\right.\]
Then either $S_1$ has a realization containing a $K_4-e$ or $S_1=(3^6)$.
\par
If $n$ is even, $S$ is a $n$-term graphical sequence, $\sigma (S)\geq 3n-2$. Let
$G$ be a   realization of $S$. Assume $d_1 \geq d_2 \geq ...\geq d_n \geq 0.$
\par
Case 1: Suppose $\sigma (S)=3n-2$. If $d_n \leq 1$, let $S'$ be the degree
sequence of $G-v_n$. Then $\sigma (S')\geq 3n-2-2=3(n-1)-1$. By induction,
$S'$ has a realization containing a $K_4-e$. Therefore $S $ has a realization
containing a $K_4-e$. Hence, we may assume that $d_n \geq
2$. Since $\sigma (S)=3n-2$, then $d_n=d_{n-1}=2$. Let $v_n$ be adjacent to
$x$ and $y$.
\par
If $x$ or $y=v_{n-1}$, let $S"$ be the degree sequence of $G-v_n-v_{n-1}$, then
$\sigma (S")=3n-2-6=3(n-2)-2$. Clearly $S"\not= (3^6)$. By induction, $S"$
has a realization containing a $K_4-e$. Hence, $S$ has a realization containing
a $K_4-e$.
\par
If $x\not=v_{n-1}$ and $y\not=v_{n-1}$, $v_{n-1}$ is adjacent to $x$ and $y$. We
now assume that $x$ is adjacent to $y$. Then $G$ contains a $K_4-e$. Hence,
we may assume that $x$ is not adjacent to $y$. Then the edge
interchange that removes the edges $xv_{n-1}$ and $yv_n$ and inserts   the
edges $xy$, $v_n v_{n-1}$ produces a realization $G'$ of $S$ containing $v_{n-1}
v_n$, and we are done as before.
\par
If $x\not=v_{n-1}$ and $y\not=v_{n-1}$, $v_{n-1}$ is not adjacent to $x$. Let $v_{n-
1}$ be adjacent to $z_1$ and $z_2$. We first consider the case $x$ is not
adjacent to $z_1$. Then the edge interchange that removes the edges
$v_{n-1}z_1$ and $v_nx$ and  inserts the edges $xz_1$ and $v_{n-1}v_n$ produces
a realization $ G'$ of $S$ containing   $v_{n-1}v_n$. We have reduced this case
to a graph $G'$ as above. Next, if $x$ is not adjacent to $z_2$. Similar to
previous case, we can prove that $S$ has a realization containing a
$K_4-e$. Finally, if  $x$ is adjacent to $z_1$ and $z_2$. We now assume that
$z_1$ is adjacent to $z_2$. Then $G$ contains a $K_4-e$. Hence, we may assume that
$z_1$ is not adjacent to $z_2$. Then the edge interchange
that removes the edges $v_{n-1}z_1$, $v_{n-1}z_2$ and $v_nx$ and inserts the
edges $v_{n-1}v_n$, $z_1z_2$ and $v_{n-1}x$ produces a realization $G'$ of  $S$
containing $v_{n-1}v_n$, and we are done as before.
\par
If $x\not=v_{n-1}$ and $y\not=v_{n-1}$, and $v_{n-1}$ is not adjacent to $y$. Similar
to  previous case, we can prove that $S$ has a realization containing a
$K_4-e$.
\par
Case 2: Suppose $\sigma (S)=3n.$
 If $d_n \leq 2$.
Let $S'$ be degree sequence of $G-v_n$, then $\sigma (S')\geq3n-4=3(n-
1)-1$. By induction, $S'$ has a realization containing a $K_4-e$.
Hence, $S$ has a realization containing a $K_4-e$.
Thus, we may assume that $d_n \geq 3.$
Then $S=(3^n).$
 If $n=6$.
Let $G_1$ be a realization of $(3^6)$. Clearly $G_1$ does not contain a $K_4-e$.
Next, if $n=4p$ $(p \geq 2).$
Then $pK_4$ is a realization of $S=(3^n)$ which contains a $K_4-e.$
 Finally, suppose that $n = 4p+2$ $ (p\geq 2).$
Then $G_1 \cup (p-1)K_4$ is a realization of $S=(3^n)$ which contains a
$K_4-e$.
\par
Case 3: Suppose $3n+2\leq \sigma (S)\leq 4n-2$.
Then $d_n \leq 3.$ Let $S'$ be a degree sequence of $G-v_n$, then $\sigma
(S')\geq 3n+2-6=3(n-1)-1$. By induction, $S'$ has a realization
containing a $K_4-e$. Hence, $S$ has a realization containing a $K_4-e$.
\par
Case 4: Suppose $\sigma (S)\geq 4n$.
If $n\geq 8$.
By [3] proposition 2 and theorem 4, $S$ has a realization containing a $K_4$.
Next, if $n=6$
and if $4n\leq \sigma (S)\leq 5n-2$.
Then $d_n\leq 4.$ Let $S'$ be a degree sequence of $G-v_n$, then $\sigma
(S')\geq 4n-8=16=3(n-1)+1.$ By induction, $S'$ has a realization
containing a $K_4-e$. Hence, $S$ has a realization containing a $K_4-e$.
Finally, Suppose that $\sigma (S)\geq  5n=30$.
Then $\sigma (S)=30$. The realization of $S$ is $K_6$ which contains $K_4-e$.
\par
If $n$ is odd, $S$ is a $n$-term graphical sequence, $\sigma (S)\geq 3n-1$.
Let $G$ be a  realization of $S$.
\par
Case 1: Suppose $\sigma (S)=3n-1$.
Then $d_n \leq 2$. Let $S'$ be degree sequence of $G-v_n$, then $\sigma
(S')\geq 3n-1-4=3(n-1)-2$. By induction, either $S'$ has
a realization containing a $K_4-e$ or $S'=(3^6)$. Therefore $S$ has a
realization containing a $K_4-e$ or $S=(4^1, 3^5,1^1)$. Clearly,
$(4^1,3^5,1^1)$ has  a realization containing a $K_4-e$ (see Appendix
Figure 1). Hence, $S$ has a  realization containing a $K_4-e$.
\par
Case 2: Suppose $3n+1\leq \sigma (S)\leq 4n-2$.
Then $d_n \leq 3$. Let $S'$ be degree sequence of $G-v_n$, then $\sigma
(S')\geq 3n+1-6=3(n-1)-2$. By induction,  either
$S'$ has a realization containing a $K_4-e$ or $S'=(3^6)$. Therefore $S$
 has a realization containing a $K_4-e$ or $S=(4^2,3^4,2^1), S=(4^3,3^4)$.
Clearly, $(4^2,3^4,2^1)$ and $(4^3,3^4)$ both have a realization containing
a $K_4-e$ (see Appendix Figure 2). Hence, $S$ has a realization containing
a $K_4-e$.
\par
Case 3: Suppose $\sigma (S)\geq 4n$. If $n\geq 9$, then
by theorem 4 of [3], $S$ has a realization containing a $K_4$.
Next, if $n=7$ and
 if $4n\leq \sigma (S)\leq 5n-1$,
then $d_n \leq 4$. Let $S'$ be a degree sequence of $G-v_n$, then $\sigma
(S') \geq 4n-8=3n-1=3(n-1)+2$. Clearly $S'\not=(3^6)$, so by induction,
$S'$ has a realization containing $K_4-e$. Thus $S$ has a realization
containing a $K_4-e$ .
Finally, Suppose  $\sigma (S)\geq 5n+1=36$. Clearly, $(6^6,0^1)$ is not
graphical. Hence $d_7\geq 1$ and by theorem 2.2 of [6], $S$ has a realization
containing a $K_4$.

\section*{Acknowledgment}
   This paper was written in the University of Science and Technology of
China as a visiting scholar. I thank Prof. Li Jiong-Sheng for his advice.
I thank the referee for many helpful comments.
\par

\section * {Appendix}

\setlength{\unitlength}{1.5mm}
\begin{picture}(60,30)
\put(15, 20){\circle*{2}}
\put(25, 20){\circle*{2}}
\put(35, 20){\circle*{2}}
\put(15, 30){\circle*{2}}
\put(25, 30){\circle*{2}}
\put(35, 30){\circle*{2}}
\put(45, 30){\circle*{2}}
\put(15, 30){\line(1, 0){30.0}}
\put(15, 20){\line(1, 0){20.0}}
\put(15, 20){\line(1, 1){10.0}}
\put(15, 30){\line(1, -1){10.0}}
\put(25, 30){\line(1, -1){10.0}}
\put(15, 20){\line(0, 1){10.0}}
\put(35, 20){\line(0, 1){10.0}}
\put(20, 15){\makebox(8, 1)[l]{$(4^1, 3^5, 1^1)$}}
\put(20, 10){\makebox(8, 1)[l]{Figure 1}}
\end{picture}
\vskip 0.1in
\setlength{\unitlength}{1.5mm}
\begin{picture}(120,40)
\put(0, 20){\circle*{2}}
\put(10, 20){\circle*{2}}
\put(20, 20){\circle*{2}}
\put(0, 30){\circle*{2}}
\put(10, 30){\circle*{2}}
\put(20, 30){\circle*{2}}
\put(30, 30){\circle*{2}}
\put(0, 30){\line(1, 0){30.0}}
\put(0, 20){\line(1, 0){20.0}}
\put(0, 20){\line(1, 1){10.0}}
\put(0, 30){\line(1, -1){10.0}}
\put(10, 30){\line(1, -1){10.0}}
\put(0, 20){\line(0, 1){10.0}}
\put(20, 20){\line(1, 1){10.0}}
\put(20, 20){\line(0, 1){10.0}}
\put(5, 10){\makebox(8, 1)[l]{$(4^2, 3^4, 2^1)$}}
\put(40, 20){\circle*{2}}
\put(50, 20){\circle*{2}}
\put(60, 20){\circle*{2}}
\put(40, 30){\circle*{2}}
\put(50, 30){\circle*{2}}
\put(60, 30){\circle*{2}}
\put(70, 30){\circle*{2}}
\put(40, 30){\line(1, 0){30.0}}
\put(40, 20){\line(1, 0){20.0}}
\put(40, 20){\line(1, 1){10.0}}
\put(40, 30){\line(1, -1){10.0}}
\put(60, 20){\line(1, 1){10.0}}
\put(50, 20){\oval(20, 10)[b]}
\put(55, 30){\oval(30, 10)[t]}
\put(40, 20){\line(0, 1){10.0}}
\put(60, 20){\line(0, 1){10.0}}
\put(45, 10){\makebox(8, 1)[l]{$(4^3, 3^4)$}}
\put(25, 5){\makebox(8, 1)[l]{Figure 2}}
\end{picture}
\end{document}